\documentstyle[11pt, amsfonts]{amsart}

\def\R{{\mathbb R}}
\def\Rad{{\mathcal R}}

\newtheorem{theorem}{Theorem}
\newtheorem{lemma}{Lemma}

\newtheorem{co}{Corollary}

\def\N{{\mathbb N}}

\def\e{\varepsilon}

\def\pf{\noindent {\bf Proof:  }}
\def\endpf{ \begin{flushright}
$ \Box $ \\
\end{flushright}}

\begin{document}

\title[Radon transform on convex sets]{Inequalities for the Radon transform on convex sets}

\author{Apostolos Giannopoulos}
\address[Apostolos Giannopoulos]{Department of Mathematics, National and Kapodistrian University of Athens,
Panepistimiopolis 157-84, Athens, Greece}
\email{apgiannop@@math.uoa.gr}

\author{Alexander  Koldobsky}
\address[Alexander  Koldobsky]{Department of Mathematics, University of Missouri, Columbia, MO 65211, USA}
\email{koldobskiya@@missouri.edu}

\author{Artem Zvavitch}
\address[Artem Zvavitch]{Department of Mathematical Sciences\\ Kent State University\\ Kent, OH USA}
\email{zvavitch@@math.kent.edu}

\thanks{Part of the work was done when the first named author visited the University of Missouri as a Miller Scholar.
The second named author is supported by the Hellenic Foundation for Research and Innovation (Project Number: 1849).
The second named author was supported in part by the U.S. National Science Foundation Grant DMS-1700036.
The third named author was supported in part by  the U.S. National Science Foundation Grant DMS-2000304 and  United States - Israel Binational Science Foundation (BSF)}
\date{}

\begin{abstract} Several years ago the authors started looking at some problems of convex geometry from a more general point of view, replacing volume by an arbitrary measure. This approach led to new general properties of the Radon transform on convex bodies including an extension of the Busemann-Petty problem and a slicing inequality for arbitrary functions.
The latter means that the sup-norm of the Radon transform of any probability density on a convex body of volume one is bounded from below by a positive constant depending only on the dimension. In this note, we prove an inequality that serves as an umbrella for these results. Let $K$ and $L$ be star bodies in $\R^n,$ let $0<k<n,$
and let $f,g$ be non-negative continuous functions on $K$ and $L$, respectively, so that $\|g\|_\infty=g(0)=1.$
Then
$$
\frac{\int_Kf}{\left(\int_L g\right)^{\frac{n-k}n}|K|^{\frac kn}}  \le \frac n{n-k} \left(d_{\rm ovr}(K,{\cal{BP}}_k^n)\right)^k
\max_{H} \frac{\int_{K\cap H} f}{\int_{L\cap H} g},
$$
where $|K|$ stands for volume of proper dimension, $C$ is an absolute constant, maximum is taken over all $(n-k)$-dimensional subspaces of $\R^n,$ and $d_{\rm ovr}(K,{\cal{BP}}_k^n)$ is the outer volume ratio distance from $K$ to the class of generalized $k$-intersection bodies in $\R^n.$
This result also implies a mean value inequality for the Radon transform, as follows. If $K$ is an origin-symmetric convex body in $\R^n,$ $0<k<n,$ and $f$ is an integrable non-negative function on $K,$ then
$$\frac{\int_Kf}{|K|}\le C^k \left(\sqrt{\frac {n}{k}\log^3\left(\frac{en}{k}\right)}\right)^k \max_H \frac{\int_{K\cap H}f}{|K\cap H|}.$$

\end{abstract}
\maketitle

\section{Introduction}\label{section:1}

The Busemann-Petty problem \cite{BP} was posed in 1956 and asks the following question. Let $K$ and $L$ be origin-symmetric convex bodies in $\R^n,$ and suppose that the $(n-1)$-dimensional volume of every central hyperplane section of $K$ is smaller than the corresponding one for $L,$ i.e. $|K\cap \xi^\bot|\le |L\cap\xi^\bot|$ for every $\xi\in S^{n-1}.$ Does it necessarily follow that the $n$-dimensional volume of $K$ is smaller than the volume of $L,$ i.e. $|K|\le |L|?$ Here $\xi^\bot=\{x\in \R^n:\ \langle x,\xi\rangle =0\}$ is the central hyperplane perpendicular to $\xi\in S^{n-1},$ and $|K|$ stands for volume of proper dimension. The answer is affirmative if the dimension $n\le 4,$ and it is negative when $n\ge 5;$ see
\cite{G3, K1} for the solution of the problem and its history.

Since the answer to the Busemann-Petty problem is negative in most dimensions, it is natural to ask
whether the inequality for volumes holds up to an absolute constant. This is known as the isomorphic
Busemann-Petty problem introduced in \cite{MP}. Does there exist an absolute constant $C$
so that for any dimension $n$ and any pair of origin-symmetric
convex bodies $K$ and $L$ in $\R^n$ satisfying
$|K\cap \xi^\bot|\le |L\cap \xi^\bot|$ for all $\xi\in S^{n-1},$
we have $ |K|\le C |L| ?$

As pointed out in \cite{MP}, the isomorphic Busemann-Petty problem is equivalent to the slicing problem of Bourgain  \cite{Bo1, Bo2}.
Does there exist an absolute constant $C$ so that for any $n\in \N$ and any
origin-symmetric convex body $K$ in $\R^n$
\begin{equation} \label{hyper}
|K|^{\frac {n-1}n} \le C \max_{\xi\in S^{n-1}} |K\cap \xi^\bot| \ ?
\end{equation}
In other words, is it true that every origin-symmetric convex body $K$ of volume one in $\R^n$ has a hyperplane section with area greater than an absolute constant, i.e. there exists $\xi\in S^{n-1}$ so that $|K\cap \xi^\bot|>c,$ where $c$ does not depend on $K$ and $n?$
The isomorphic Busemann-Petty problem and the slicing problem are still open. Klartag \cite{Kla} proved that $C\le O(n^{1/4})$, removing a logarithmic term from an earlier estimate of Bourgain \cite{Bo3}.

Several years ago the authors decided to look at these problems from a more general point of view, replacing volume by an arbitrary measure. An extension of the slicing problem was considered in \cite{K2, K3, K4, CGL, KK, KLi}. It was proved in \cite{K4} that for any $n\in \N,$ any star body $K$ in $\R^n$ and any non-negative continuous function $f$ on $K,$
\begin{equation}\label{main-problem}
\int_Kf\ \le\ 2\ d_{\rm {ovr}}(K,{\cal{I}}_n)\  |K|^{1/n}\ \max_{\xi \in S^{n-1}} \int_{K\cap \xi^\bot} f.
\end{equation}
Here
\begin{equation}\label{eq:dvr}d_{{\rm {ovr}}}(K,   \Omega) = \inf \left\{ \left( \frac {|D|}{|K|}\right)^{1/n}:\ K\subset D,\ D\in  \Omega \right\}\end{equation}
is the outer volume ratio distance from $K$ to a fixed  class $\Omega$ of bodies in $\R^n$, in the particular case of equation (\ref{main-problem}) 
that we consider, $\Omega={\cal{I}}_n$ is the class of intersection bodies (see definition below).

It is interesting to note that we do not need to require the function $f$ to be even, as well as the body $K$ to be symmetric, and additional assumptions are only needed to estimate $d_{\rm {ovr}}(K,{\cal{I}}_n)$.  Since the class of intersection bodies includes ellipsoids, by John's theorem \cite{J}, if $K$ is origin-symmetric and convex, then $d_{{\rm {ovr}}}(K,{\cal{I}}_n)\le \sqrt{n}.$
Thus, there exists a constant $s_n\ge \frac1{2\sqrt{n}}$ so that for any origin-symmetric convex body $K$ of volume 1 in $\R^n$ and any integrable non-negative function $f$ on $K$ with $\int_Kf=1,$ there exists a direction $\xi\in S^{n-1}$ for which $\int_{K\cap \xi^\bot} f \ge s_n.$ In other words, the sup-norm of the Radon transform of any probability density on a convex body of volume one is bounded from below by a positive constant depending only on the dimension. Note that this estimate was extended later to the case of non-symmetric bodies in \cite{CGL}. An extension to the derivatives of the Radon transform was obtained in \cite{GK}.

On the other hand, it was proved in \cite{KK} that there exists an origin-symmetric convex body $M$
in $\R^n$ and a probability density $f$ on $M$ so that
$$\int_{M\cap H}f \le C\frac{\sqrt{\log\log n}}{\sqrt{n}} |M|^{-1/n},$$
for every affine hyperplane $H$ in $\R^n,$ where $C$ is an absolute constant. The logarithmic term was
later removed in \cite{KLi}, so $s_n\le  C/\sqrt{n}.$ Finally, $ \frac{c_1}{\sqrt{n}} \le s_n \le \frac{c_2}{\sqrt{n}},$   where $c_1,c_2>0$ are absolute constants.

A lower dimensional version of inequality (\ref{main-problem}) was also proved in \cite{K4}. If $K$ is a star body in $\R^n,$ $f$ is a continuous non-negative function on $K$, and $1\le k <n,$ then
\begin{equation}\label{main-problem1}
\int_Kf\ \le\ C^k \ (d_{\rm {ovr}}(K,{\cal{BP}}_k^n))^k\  |K|^{k/n}\ \max_{H\in Gr_{n-k}} \int_{K\cap H} f,
\end{equation}
where $C$ is an absolute constant, and ${\cal{BP}}_k^n$ is the class of generalized $k$-intersection bodies (see definition below). This implies that there exists a constant $c_{n,k}>0$ such that for any convex body $K$ in $\R^n$ and any probability density $f$ on $K,$ there exists an $(n-k)$-dimensional subspace $H$ in $\R^n$ so that
$\int_{K\cap H}f \ge c_{n,k}.$ Moreover, applying inequality (\ref{kpzy}) we get
$$(c_{n,k})^{1/k}\ge \frac{c\sqrt{k}}{\sqrt{n\log^3(\frac{en}k)}},$$
where $c>0$ is an absolute constant. It is an open problem whether it is possible to remove the logarithmic term
in this estimate or to find the exact upper estimate for $c_{n,k}.$

An extension of the Busemann-Petty problem to arbitrary functions was found in \cite{Zv1, Zv2}. Suppose that $f$ is an even continuous strictly positive function on $\R^n,$ and $K$ and $L$ are origin-symmetric convex bodies in $\R^n$ so that
\begin{equation}\label{bp-f}
\int_{K\cap \xi^\bot} f \le \int_{L\cap \xi^\bot} f, \qquad \forall \xi\in S^{n-1}.
\end{equation}
Does it necessarily follow that $\int_K f\le \int_Lf?$ The answer is the same as for volume:
affirmative if $n\le 4$ and negative if $n\ge 5.$

An isomorphic version was proved in \cite{KZ}.
For every dimension $n,$ inequalities (\ref{bp-f}) imply $\int_K f\le \sqrt{n}\ \int_Lf.$
In fact, it was proved in \cite{KZ} that inequalities (\ref{bp-f}) imply
$\int_K f\le   d_{BM}(K,{\cal{I}}_n) \int_Lf,$ where
$$d_{BM}(K,{\cal{I}}_n)=\inf \{a>0: \exists D \in {\cal{I}}_n: D \subset K \subset aD \}$$
is the Banach-Mazur distance from $K$ to the class of intersection bodies.
Now if $K$ is origin-symmetric and convex, by John's theorem, $d_{BM}(K,{\cal{I}}_n)\le \sqrt{n},$ so the $\sqrt{n}$ estimate follows. It is not known whether the $\sqrt{n}$ estimate is optimal in this problem.
Another version of the isomorphic Busemann-Petty problem was proved in \cite{KPZv}; see Section \ref{section:5}.


In this article we prove a theorem that generalizes these results.

\begin{theorem} \label{quotient-main} {} Let $K$ and $L$ be star bodies in $\R^n,$ let $0<k<n,$
and let $f,g$ be non-negative continuous functions on $K$ and $L$, respectively, so that $\|g\|_\infty=g(0)=1.$
Then
\begin{equation}\label{quotient}
\frac{\int_Kf}{\left(\int_L g\right)^{\frac{n-k}n}|K|^{\frac kn}}  \le \frac n{n-k} \left(d_{\rm ovr}(K,{\cal{BP}}_k^n)\right)^k
\max_{H\in Gr_{n-k}} \frac{\int_{K\cap H} f}{\int_{L\cap H} g}.
\end{equation}
\end{theorem}

In Section 5, we deduce various properties of the Radon transform from the latter inequality.
For example, if we put $K=L$ and $g\equiv 1$ in (\ref{quotient}), we get what we call {the mean value inequality for the Radon transform}:
\begin{equation}\label{mean}
\frac{\int_Kf}{|K|}\le \frac n{n-k} \left(d_{\rm ovr}(K,{\cal{BP}}_k^n)\right)^k \max_{H\in Gr_{n-k}} \frac{\int_{K\cap H}f}{|K\cap H|}.
\end{equation}
The distance  $d_{\rm ovr}(K,{\cal{BP}}_k^n)$ does not depend on $n$ for some special classes of bodies $K.$ This distance is equal to 1 for intersection bodies $K,$ because ${\cal{I}}_n\subset {\cal{BP}}_k^n$ for all $k$ \cite{GZ, M2}. It was proved in \cite{K4} that for unconditional convex bodies $K,$ $d_{{\rm {ovr}}}(K,{\cal{I}}_n)\le e.$  Also, this distance is bounded by an absolute constant for the polar bodies of convex bodies with bounded volume ratio \cite{K4}.
For the unit balls of $n$-dimensional subspaces of $L_p,\ p>2,$ the distance is less than $c\sqrt{p},$ where $c>0$ is an absolute constant, as shown in \cite{M1, KP}. Note that the unit balls of subspaces of $L_p$ with $0<p\le 2$ are intersection bodies \cite{K5}. In general, it was proved in \cite{KPZy} that for any origin-symmetric convex body $K$
\begin{equation}\label{kpzy}
d_{{\rm {ovr}}}(K,{\cal BP}_k^n)\le C\sqrt{\frac nk}\log^{\frac 32}\left(\frac {en}k\right),
\end{equation}
where $C$ is an absolute constant. Thus, if $K$ is origin symmetric convex, the mean value inequality holds with a constant depending only on $n,k:$
$$\frac{\int_Kf}{|K|}\le C^k \left(\sqrt{\frac{n}{k}\log^3\left(\frac{en}{k}\right )}\right)^k \max_H \frac{\int_{K\cap H}f}{|K\cap H|}.$$

In Section~\ref{section:3} we describe an alternative approach to Theorem~\ref{quotient-main} which is based on Blaschke-Petkantchin
formulas and affine isoperimetric inequalities. This approach was initiated in \cite{CGL} and leads to a version of \eqref{quotient}
which is valid for more general pairs of sets. Below, for any bounded Borel set $K$ in ${\mathbb R}^n$ we denote by ${\rm ovr}(K)$
the outer volume ratio ${\rm ovr}(K)=d_{\rm ovr}(K, L_2)=\inf\limits_{{\cal E}}\big(|{\cal E}|/|K|\big)^{1/n}$, where the infimum is over all origin
symmetric ellipsoids ${\cal E}$ in ${\mathbb R}^n$ with $K\subseteq {\cal E}$.

\begin{theorem}\label{th:arb-ovr}Let $K$ and $L$ be two bounded Borel sets in ${\mathbb R}^n$. Let $f$ and $g$ be two bounded
non-negative measurable functions on $K$ and $L$, respectively, and assume that $\|g\|_1>0$ and $\|g\|_{\infty }=1$. For every $1\leq k\leq n-1$
we have that
\begin{equation}\frac{\int_Kf}{\left(\int_L g\right)^{\frac{n-k}n}|K|^{\frac kn}}  \le \left(C\cdot {\rm ovr}(K)\right)^k
\max_{H\in Gr_{n-k}} \frac{\int_{K\cap H} f}{\int_{L\cap H} g},\end{equation}
where $C>0$ is an absolute constant.
\end{theorem}

It should be noted that even for origin-symmetric convex bodies $K$ in ${\mathbb R}^n$ the outer volume ratio
${\rm ovr}(K)$ can be as large as $\sqrt{n}$. This is a disadvantage of Theorem~\ref{th:arb-ovr} which
does not provide estimates depending on $k$ in contrast to Theorem~\ref{quotient-main}. Regarding this
comparison, we mention that besides Blaschke-Petkantchin formulas and affine isoperimetric inequalities, the proof of Theorem~\ref{th:arb-ovr} exploits
a well-known result of Bar\'{a}ny and F\"{u}redi \cite{Barany-Furedi-1988} which may be stated as follows:
if ${\cal E}$ is an ellipsoid in ${\mathbb R}^m$, $s\geq m+1$ and $w_1,\ldots ,w_m\in {\cal E}$, then
\begin{equation*}\left (\frac{|{\rm conv}(w_1,\ldots ,w_m)|}{|{\cal E}|}\right)^{1/m}\leq C \sqrt{\log
(1+s/m)/m}\end{equation*}
where $C>0$ is an absolute constant. It would be interesting to obtain an optimal estimate for the corresponding
result when $w_1,\ldots ,w_m$ are chosen from a body $L\in {\cal{BP}}_k^n$.

In Section~\ref{section:4} we obtain a generalization of the isomorphic version of the Shephard problem due to Ball \cite{Ba2,Ba3}. Ball has proved that if $K$ and $L$ are origin-symmetric convex bodies in $\R^n$ such that
$$|K\vert\xi^\bot|\le |L\vert\xi^\bot|,\qquad \forall \xi\in S^{n-1},$$
then
$$|K|\le d_{\rm vr}(L,\Pi_n) |L|.$$
In this statement, $K\vert\xi^\bot$ denotes the orthogonal projection of $K$ in the direction of $\xi $ and
$d_{{\rm {vr}}}(L,\Pi_n)$, defined in \eqref{eq:dvr}, is the  volume ratio distance from $L$ to the class $\Pi_n$ 
of projection bodies. Replacing $\Pi_n$ by the class $\Pi_{p,n}$ of $p$-projection bodies (see Section~\ref{section:4} 
for background information) we obtain the following:

\begin{theorem} \label{main-proj} Fix $p\ge 1$ and let $K$ and $L$ be
convex bodies in $\R^n,$ then
$$\left(\frac{|K|}{|L|}\right)^{\frac{n-p}{pn}}  \le \ d_{\rm vr} (L, \Pi_{p,n})\max_{\xi\in S^{n-1}} \frac{h_{\Pi_p K} (\xi)}{h_{\Pi_p L} (\xi)}.$$
\end{theorem}

Throughout the paper, we write $a\simeq b$ if $c_1b\le a \le c_2b,$ where $c_1,c_2>0$ are absolute constants.

\section{Proof of Theorem \ref{quotient-main}}\label{section:2}

We need several definitions and facts.
A closed bounded set $K$ in $\R^n$ is called a {\it star body}  if
every straight line passing through the origin crosses the boundary of $K$
at exactly two points different from the origin, the origin is an interior point of $K,$
and the {\it Minkowski functional}
of $K$ defined by
$$\|x\|_K = \min\{a\ge 0:\ x\in aK\}$$
is a continuous function on $\R^n.$
We use the polar formula for the volume $|K|$ of a star body $K:$
\begin{equation}\label{polar-vol}
|K|=\frac 1n \int_{S^{n-1}} \|\theta\|_K^{-n} d\theta.
\end{equation}

If $f$ is an integrable function on $K$, then
\begin{equation}\label{polar-meas}
\int_K f = \int_{S^{n-1}} \left(\int_0^{\|\theta\|_K^{-1}} r^{n-1}f(r\theta) dr \right) d\theta.
\end{equation}

Let $Gr_{n-k}$ be the Grassmanian of $(n-k)$-dimensional subspaces of $\R^n.$
For $1\le k \le n-1,$  the {\it $(n-k)$-dimensional spherical Radon transform}
$\Rad_{n-k}:C(S^{n-1})\to C(Gr_{n-k})$ is a linear operator defined by
$$\Rad_{n-k}g (H)=\int_{S^{n-1}\cap H} g(x)\ dx,\quad \forall  H\in Gr_{n-k}$$
for every function $g\in C(S^{n-1}).$
\smallbreak
For every $H\in Gr_{n-k},$ the $(n-k)$-dimensional volume of the section of a star body $K$ by $H$
can be written as
\begin{equation}\label{vol-sect}
|K\cap H| = \frac1{n-k} \Rad_{n-k}(\|\cdot\|_K^{-n+k})(H).
\end{equation}
More generally, for an integrable function $f$ and any $H\in Gr_{n-k}$,
\begin{equation}\label{meas-sect}
\int_{K\cap H} f = \Rad_{n-k}\left(\int_0^{\|\cdot\|_K^{-1}} r^{n-k-1}f(r\ \cdot)\ dr \right)(H).
\end{equation}

The class of intersection bodies ${\cal{I}}_n$ was introduced by Lutwak \cite{L2}. We consider 
a generalization of this concept due to Zhang \cite{Z2}. We say that an origin symmetric star body $D$ 
in $\R^n$ is a {\it generalized $k$-intersection body},
and write $D\in {\cal{BP}}_k^n,$  if there exists a finite Borel non-negative measure $\nu_D$
on $Gr_{n-k}$ so that for every $g\in C(S^{n-1})$
\begin{equation}\label{genint}
\int_{S^{n-1}} \|x\|_D^{-k} g(x)\ dx=\int_{Gr_{n-k}} R_{n-k}g(H)\ d\nu_D(H).
\end{equation}
When $k=1$ we get the original Lutwak's class of {\it intersection bodies} ${\cal{BP}}_1^n={\cal{I}}_n$.
\bigbreak

\noindent{\bf Proof of Theorem \ref{quotient-main}.} For a small $\delta>0,$ let $D\in {\cal{BP}}_k^n$ be a body such that $K\subset D$ and
\begin{equation}\label{sect11}
|D|^{\frac 1n}\le (1+\delta)\ d_{\rm ovr}(K,{\cal{BP}}_k^n)\ |K|^{\frac 1n},
\end{equation}
and let $\nu_D$ be the measure on $Gr_{n-k}$ corresponding to $D$ by the definition (\ref{genint}).

Let $\e$ be such that
$$\int_{K\cap H} f\le \e \int_{L\cap H} g,\qquad \forall H\in Gr_{n-k}.$$

By (\ref{meas-sect}), we have
$$\Rad_{n-k}\left(\int_0^{\|\cdot\|_K^{-1}} r^{n-k-1}f(r\ \cdot)\ dr \right)(H) \le
\e\ \Rad_{n-k}\left(\int_0^{\|\cdot\|_L^{-1}} r^{n-k-1}g(r\ \cdot)\ dr \right)(H)$$
for every $H\in Gr_{n-k}.$
Integrating both sides of the latter inequality with respect to $\nu_D$ and using the definition (\ref{genint}), we get
\begin{align}\label{integration}
&\int_{S^{n-1}} \|x\|_D^{-k} \left(\int_0^{\|x\|_K^{-1}} r^{n-k-1}f(rx)\ dr \right)dx\\
\nonumber &\hspace*{1cm}\le \e \int_{S^{n-1}} \|x\|_D^{-k} \left(\int_0^{\|x\|_L^{-1}} r^{n-k-1}g(rx)\ dr \right)dx,
\end{align}
which is equivalent to
\begin{equation}\label{sect12}
\int_K \|x\|_D^{-k}f(x)dx \le \e \int_L \|x\|_D^{-k}g(x)dx.
\end{equation}
Since $K\subset D,$ we have $1\ge \|x\|_K\ge \|x\|_D$ for every $x\in K.$ Therefore,
$$\int_K \|x\|_D^{-k}f(x) dx \ge \int_K \|x\|_K^{-k}f(x) dx \ge \int_K f.$$
On the other hand, by
\cite[Lemma 2.1]{MP},
$$
\left(\frac{\int_{L}\|x\|_D^{-k} {g}(x) dx}{ \int_{D}\|x\|_D^{-k} dx} \right)^{1/(n-k)}  \le \left(\frac{\int_{L}{g}(x) dx \large}{  \int_{D} dx} \right)^{1/n}.
$$
Since $\int_D\|x\|_D^{-k} dx =\frac{n}{n-k} |D|,$ we can estimate the right-hand side of
(\ref{sect12}) by
$$\int_L \|x\|_D^{-k} g(x) dx \le\e \frac n{n-k}  \left(\int_Lg\right)^{\frac {n-k}n} |D|^{\frac kn}.$$
Applying (\ref{sect11}) and sending $\delta$ to zero, we see that the latter inequality in conjunction with (\ref{sect12})
implies
$$\int_K f \le \e\ \frac n{n-k} \left(d_{\rm ovr}(K,{\cal{BP}}_k^n)\right)^k |K|^{\frac kn}.$$
Now put
$\e = \max\limits_{H\in Gr_{n-k}} \frac{\int_{K\cap H} f}{\int_{L\cap H} g}.$ \endpf

If $f\equiv 1,\ g\equiv 1,$ we can get a slightly sharper inequality than what Theorem \ref{quotient-main} gives in this case.

\begin{theorem} \label{quotient-holder}Let $K,L$ be star bodies in $\R^n$ and $0<k<n,$ then
$$\left(\frac{|K|}{|L|}\right)^{\frac {n-k}n}\le \left(d_{\rm ovr}(K,{\cal{BP}}_k^n)\right)^k \max_{H\in Gr_{n-k}} \frac{|K\cap H|}{|L\cap H|}.$$
\end{theorem}

\pf Let $\e$ be such that $|K\cap H|\le\e\ |L\cap H|$ for all $H\in Gr_{n-k},$ and let $D$ be as in the proof of Theorem
\ref{quotient}. By (\ref{vol-sect}), for all $H$
$$\Rad_{n-k}(\|\cdot\|_K^{-n+k})(H)\le \e\ \Rad_{n-k}(\|\cdot\|_L^{-n+k})(H).$$
Integrating both sides with respect to $\nu_D$ and using the definition (\ref{genint}) we get
$$\int_{S^{n-1}} \|x\|_D^{-k}\|x\|_K^{-n+k} dx \le\e  \int_{S^{n-1}} \|x\|_D^{-k}\|x\|_L^{-n+k} dx.$$
Since $K\subset D,$ we have $1\ge \|x\|_K\ge \|x\|_D,$ and by (\ref{polar-vol}) the left-hand side is greater than $n|K|.$
Using this and H\"older's inequality,
\begin{align*}n|K| &\le\e\ \int_{S^{n-1}} \|x\|_D^{-k}\|x\|_L^{-n+k} dx\le\e\ \left(\int_{S^{n-1}} \|x\|_D^{-n}dx \right)^{\frac kn}
\left(\int_{S^{n-1}} \|x\|_L^{-n}dx \right)^{\frac {n-k}n}\\
&=\e\  n |D|^{\frac kn} |L|^{\frac {n-k}n}\le \e\ n\ (1+\delta)^k\ \left(d_{\rm ovr}(K,{\cal{BP}}_k^n)\right)^k\ |K|^{\frac kn}|L|^{\frac {n-k}n}.
\end{align*}
Sending $\delta$ to zero and setting
$$\e=\max_{H\in Gr_{n-k}} \frac{|K\cap H|}{|L\cap H|},$$
we get the result. \endpf

\section{Proof of Theorem~\ref{th:arb-ovr}}\label{section:3}

Our first tool will be a Blaschke-Petkantchin formula (see \cite[Chapter 7.2]{Schneider-Weil-book} and \cite[Lemma 5.1]{Gardner-2007}).

\begin{lemma}[Blaschke-Petkantschin]\label{lem:blaschke}Let $1\leq s\leq n-1$. There exists a constant $p(n,s)>0$ such that, for every non-negative
bounded Borel measurable function $F:({\mathbb R}^n)^s\to {\mathbb R}$,
\begin{align}\label{eq:tools-1}&\int_{{\mathbb R}^n}\cdots \int_{{\mathbb R}^n}F(x_1,\ldots ,x_s)dx_s\cdots dx_1\\
\nonumber &\hspace*{1cm} =p(n,s)\int_{G_{n,s}}\int_H\cdots \int_Ff(x_1,\ldots ,x_s)\,|{\rm conv}(0,x_1,\ldots ,x_s)|^{n-s}\\
\nonumber &\hspace*{3cm}dx_s\cdots dx_1\,d\nu_{n,s}(H),
\end{align}
where $\nu_{n,s}$ is the Haar probability measure on $Gr_s$. The exact value of the constant $p(n,s)$ is
\begin{equation}\label{eq:tools-2}p(n,s)=(s!)^{n-s}\frac{(n\omega_n)\cdots ((n-s+1)\omega_{n-s+1})}{(s\omega_s)\cdots (2\omega_2)\omega_1},\end{equation}
where $\omega_n$ is the volume of the unit Euclidean ball $B_2^n.$
\end{lemma}

We shall also use the next inequality, proved independently by Busemann and Straus \cite{Busemann-Straus-1960}, and  Grinberg \cite{Grinberg-1990}.

\begin{lemma}[Busemann-Straus, Grinberg]\label{lem:grinberg}Let $K$ be a bounded Borel set of volume $1$ in ${\mathbb R}^n$.
For any $1\leq k\leq n-1$ and $T\in SL(n)$ we have
\begin{equation}\label{eq:grinberg}\int_{Gr_{n-k}}|K\cap H|^n
d\nu_{n,n-k}(H)\leq \int_{Gr_{n-k}}|\overline{B}_2^n\cap H|^n
d\nu_{n,n-k}(H),\end{equation}
where $\overline{B}_2^n$ is the Euclidean ball of volume $1$.
\end{lemma}

Our next tool will be a  theorem of Dann, Paouris and Pivovarov from \cite{Dann-Paouris-Pivovarov-2015}; 
the proof of this fact combines Blaschke-Petkantschin formulas with rearrangement inequalities.

\begin{lemma}[Dann-Paouris-Pivovarov]\label{lem:dann-1}Let $g$ be a non-negative,
bounded integrable function on ${\mathbb R}^n$ with $\|g\|_1>0$. For every $1\leq k\leq n-1$ we have
\begin{equation}\label{eq:dann-1}\int_{Gr_{n-k}}\frac{1}{\|g|_H\|_{\infty }^k}\left (\int_Hg(x)dx\right )^nd\nu_{n,n-k}(H)
\leq \gamma_{n,k}^{-n}\left (\int_{{\mathbb R}^n}g(x)dx\right )^{n-k},\end{equation}
where $\gamma_{n,k}=\omega_n^{\frac{n-k}{n}}/\omega_{n-k}$.
\end{lemma}

It is checked in \cite{CGL} that for every $1\leq k\leq n-1$ one has
\begin{equation}\label{eq:tools-13}e^{-k/2}<\gamma_{n,k}<1\quad\hbox{and}\quad [\gamma_{n,k}^{-n}p(n,n-k)]^{\frac{1}{k(n-k)}}\simeq \sqrt{n-k}.\end{equation}
Finally, we need a well-known theorem of B\'{a}r\'{a}ny and F\"{u}redi \cite{Barany-Furedi-1988}: if $s\geq m+1$ and $w_j\in {\mathbb R}^m$ satisfy $\|w_j\|_2\leq 1$
for $j=1,\ldots ,s$, then
\begin{equation*}|{\rm conv}(w_1,\ldots ,w_s)|^{1/m}\leq C \frac{\sqrt {\log
(1+s/m)}}{m}.\end{equation*}
Equivalently, this says that if $w_j\in B_2^m$, $1\leq j\leq s$, then the volume radius of their convex hull is bounded
by $C \sqrt {\log(1+s/m)/m}$. By affine invariance we obtain:

\begin{lemma}\label{lem:barany-furedi}There exists an absolute constant $C>0$ such that if ${\cal E}$ is an ellipsoid in ${\mathbb R}^m$, $s\geq m+1$ and $w_1,\ldots ,w_m\in {\cal E}$, then
\begin{equation*}\left (\frac{|{\rm conv}(w_1,\ldots ,w_m)|}{|{\cal E}|}\right)^{1/m}\leq C \sqrt{\log
(1+s/m)/m}.\end{equation*}
\end{lemma}

\noindent {\bf Proof of Theorem~\ref{th:arb-ovr}.}
Let ${\cal E}$ be a centered ellipsoid such that $K\subseteq {\cal E}$ and
$${\rm ovr}(K)=(|{\cal E}|/|K|)^{1/n}.$$
We shall use the next consequence of Lemma~\ref{lem:barany-furedi}: if $F\in Gr_{n-k}$ and $x_1,\ldots ,x_{n-k}\in K\cap H\subseteq {\cal E}\cap H$
then ${\rm conv}(0,x_1,\ldots ,x_{n-k})\subseteq {\cal E}\cap H$, and since ${\cal E}\cap H$ is an $(n-k)$-dimensional
centered ellipsoid we must have
\begin{align}\label{eq:ovr-convex-hull}|{\rm conv}(0,x_1,\ldots ,x_{n-k})| &\leq\left( C_1\frac{\sqrt{\log(1+(n-k+1)/(n-k)}}{\sqrt{n-k}}\right)^k
|{\cal E}\cap H|\\
\nonumber &\leq \left( \frac{C_1}{\sqrt{n-k}}\right)^k |{\cal E}\cap H|.\end{align}
Applying Lemma~\ref{lem:blaschke} for the function $F(x_1,\ldots ,x_{n-k})=\prod_{i=1}^{n-k}f(x_i){\bf 1}_{K}(x_i)$, with $s=n-k$, we get
\begin{align*}\left (\int_Kf(x)\,dx\right)^{n-k} &=\int_{{\mathbb R}^n}\cdots \int_{{\mathbb R}^n}F(x_1,\ldots ,x_{n-k})dx_{n-k}\cdots dx_1\\
\nonumber &=p(n,n-k)\int_{Gr_{n-k}}\int_{K\cap H}\cdots\int_{K\cap H}g(x_1)\cdots g(x_{n-k})\\
\nonumber &\hspace*{1cm}\times \,|{\rm conv}(0,x_1,\ldots ,x_{n-k})|^{k}dx_{n-k}\cdots dx_1\,d\nu_{n,n-k}(H).
\end{align*}
Let $M:=\max_{H\in Gr_{n-k}} \frac{\int_{K\cap H} f}{\int_{L\cap H} g}$. Then, \eqref{eq:ovr-convex-hull} shows that
\begin{align*}&\left (\int_Kf(x)\,dx\right)^{n-k}\leq p(n,n-k)\left (\frac{C_1}{\sqrt{n-k}}\right)^{k(n-k)}\int_{Gr_{n-k}}\int_{K\cap H}\cdots\int_{K\cap H}|{\cal E}\cap H|^{k}\\
\nonumber &\times\,f(x_1)\cdots f(x_{n-k})\,dx_{n-k}\cdots dx_1\,d\nu_{n,n-k}(H)\\
\nonumber &= p(n,n-k)\left (\frac{C_1}{\sqrt{n-k}}\right)^{k(n-k)}\int_{Gr_{n-k}}|{\cal E}\cap H|^k\Big(\int_{K\cap H}f(x)\,dx\Big)^{n-k}\,d\nu_{n,n-k}(H)\\
\nonumber &\leq p(n,n-k)\left (\frac{C_1}{\sqrt{n-k}}\right)^{k(n-k)}M^{n-k}\int_{Gr_{n-k}}|{\cal E}\cap H|^k\\
&\times \Big(\int_{L\cap H}g(x)\,dx\Big)^{n-k}\,d\nu_{n,n-k}(H).
\end{align*}
Now, by H\"{o}lder's inequality and Grinberg's inequality (\ref{eq:grinberg})  we get
\begin{align*}
&\int_{Gr_{n-k}}|{\cal E}\cap H|^k\left(\int_{L\cap H}g(x)\,dx\right)^{n-k}\,d\nu_{n,n-k}(H)\\
&\leq\left (\int_{Gr_{n-k}}|{\cal E}\cap H|^nd\nu_{n,n-k}(H)\right )^{\frac{k}{n}}
\left (\int_{Gr_{n-k}}\left(\int_{L\cap H}g(x)\,dx\right)^nd\nu_{n,n-k}(H)\right)^{\frac{n-k}{n}}\\
&\leq \gamma_{n,k}^{-n}|{\cal E}|^{\frac{k(n-k)}{n}}\left (\int_{Gr_{n-k}}\left(\int_{L\cap H}g(x)\,dx\right)^nd\nu_{n,n-k}(H)\right)^{\frac{n-k}{n}}\\
&=\gamma_{n,k}^{-n}|K|^{\frac{k(n-k)}{n}}{\rm ovr}(K)^{k(n-k)}\left (\int_{Gr_{n-k}}\left(\int_{L\cap H}g(x)\,dx\right)^nd\nu_{n,n-k}(H)\right)^{\frac{n-k}{n}}.
\end{align*}
Finally, since $\|g|_H\|_{\infty }\leq \|g\|_{\infty }=1$ for all $H\in Gr_{n-k}$, we may apply \eqref{eq:dann-1} to get
\begin{equation}\label{eq:dann-2}\int_{Gr_{n-k}}\left (\int_{L\cap H}g(x)dx\right )^nd\nu_{n,n-k}(H)
\leq \gamma_{n,k}^{-n}\left (\int_Lg(x)dx\right )^{n-k}.\end{equation}
Combining the above we get
\begin{align}\label{eq:final-ovr}\left (\int_Kf(x)\,dx\right)^{n-k}&\leq \left [\gamma_{n,k}^{-n}p(n,n-k)(C_1/\sqrt{n-k})^{k(n-k)}\right]\\
\nonumber &\hspace*{1cm}\times M^{n-k}|K|^{\frac{k(n-k)}{n}}{\rm ovr}(K)^{k(n-k)}\left (\int_Lg(x)dx\right )^{\frac{(n-k)^2}{n}}.\end{align}
Note that, by \eqref{eq:tools-13},
\begin{equation*}\left [\gamma_{n,k}^{-n}p(n,n-k)\right]^{\frac{1}{n-k}}\left (\frac{C}{\sqrt{n-k}}\right )^k\leq C^k
\end{equation*}
for some absolute constant $C>0$. Then, the result follows from \eqref{eq:final-ovr}. \hfill $\Box $

\bigskip

\section{Proof of Theorem~\ref{main-proj}}\label{section:4}

The proof of Theorem~\ref{main-proj} requires several additional definitions and facts from convex geometry. We refer the reader
to \cite{S2} for details.

The support function of a convex body $K$ in $\R^n$ is defined by
$$
h_K(x) = \max_{\xi \in K} \langle x,\xi\rangle,\quad x\in \R^n.
$$
If $K$ is origin-symmetric, then $h_K$ is a norm on $\R^n.$ One of the crucial properties of
the support function is its relation to the Minkowski sum of convex bodies:
\begin{equation}\label{supMin}
h_{K+L}(x)=h_K(x)+h_L(x).
\end{equation}
The {surface area measure} $S(K, \cdot)$ of a convex body $K$ in
$\R^n$ is defined as follows: for every Borel set $E \subset S^{n-1},$
$S(K,E)$ is equal to Lebesgue measure of the part of the boundary of $K$
where normal vectors belong to $E.$ The volume of a convex body can be expressed in terms of its support function and
surface area measure:
\begin{equation}\label{vol-proj}
|K| = \frac 1n \int_{S^{n-1}}h_K(x) dS(K, x).
\end{equation}
If $K$ and $L$ are two convex bodies in $\R^n,$ the {\it mixed volume} $V_1(K,L)$ is equal to
\begin{equation}\label{mixv}
V_1(K,L)= \frac{1}{n} \lim_{\e\to +0}\frac{|K+\epsilon L|- |K|}{\e}.
\end{equation}
We shall use the first Minkowski inequality: for any pair of convex bodies $K,L$ in $\R^n,$
\begin{equation} \label{mink}
V_1(K,L) \ge |K|^{\frac{n-1}n} |L|^{1/n}.
\end{equation}
The mixed volume $V_1(K,L)$ can also be expressed in terms of the support function and surface area measure:
\begin{equation}\label{mixed}
V_1(K,L) = \frac 1n \int_{S^{n-1}}h_L(x) d S(K,x).
\end{equation}
For a convex body $K$ in $\R^n$ and $\xi\in S^{n-1},$ denote by $K\vert \xi^\bot$ the orthogonal projection of $K$ to the central
hyperplane $\xi^\bot.$ The Cauchy formula states that
\begin{equation}\label{cauchy}
|K\vert \xi^\bot|=\frac 12 \int_{S^{n-1}} |\langle x,\xi\rangle|\ dS(K,x).
\end{equation}
Let $K$ be a convex body in $\R^n.$ The projection body $\Pi K$ of $K$ is defined as an origin-symmetric convex
body in $\R^n$ whose support function in every direction is equal to the volume of the orthogonal projection of $K$ to this direction:
for every $\theta\in S^{n-1},$
\begin{equation} \label{def:proj}
h_{\Pi K}(\xi) = |K\vert \xi^{\perp}|.
\end{equation}
We denote by $\Pi_n$ the class of projection bodies of convex bodies and if $D\in\Pi_n$ we simply say
that $D$ is a projection body. By Cauchy's formula (\ref{cauchy}), for every projection body $D$
there exists a finite measure $\nu_D$ on $S^{n-1}$ such that
\begin{equation}\label{proj-body}
h_D(x)=\int_{S^{n-1}} |\langle x,\xi\rangle|\ d\nu_D(\xi),\qquad \forall x\in S^{n-1}.
\end{equation}
Let ${\cal K}_0$ denote the class of convex bodies containing the origin in their interior. Firey \cite{F} extended the concept of Minkowski sum (\ref{supMin}), and introduced for each real $p \ge 1$, a new linear combination
of convex bodies, the so-called $p$-sum:
$$
h^p_{\alpha  K +_{p} \beta  L} (x)
=\alpha h^p_{K}(x) +\beta h^p_{L}(x).
$$
Here $K,L \in {\cal K}_0$ and $\alpha, \beta $ are positive real numbers.
In a series  of papers Lutwak \cite{Lu2, Lu3} showed that the  Firey sums lead to a
Brunn-Minkowski theory for each $p \ge 1$. Extending the classical definition of the mixed volume (\ref{mixv})
Lutwak introduced the notion of $p$-mixed volume, $V_p(K,L)$, $p \ge 1$ as
$$
\,V_p(K,L)= \frac{p}{n} \lim\limits_{\e \to 0}\frac{V(K+_p\e \, L)-V(K)}{\e},
$$
for all $K, L \in  {\cal K}_0$. Lutwak proved that for each $K \in {\cal K}_0$, there exists a positive Borel
measure $S_p(K, \cdot)$ on $S^{n-1}$ so that
$$
V_p(K,L)=\frac{1}{n} \int\limits_{S^{n-1}} h^p_L(u) dS_p(K, u)
$$
for all $L \in {\cal K}_0$. It turns out that the measure $S_p(K, \cdot)$ is absolutely continuous
with respect to $S(K, \cdot)$, with Radon-Nikodym derivative
$$
\frac{dS_p(K, \cdot)}{dS(K, \cdot)}\,=\,h(K, \cdot)^{1-p}.
$$
Lutwak \cite{Lu2} generalized  the Brunn-Minkowski inequality to the case of $p$-mixed volumes as follows:
\begin{equation}\label{pmin}
V_p(K, L)^n \,\ge \,|K|^{n-p} |L|^{p}, \qquad p>1.
\end{equation}
We will also use the concept of a $p$-projection body, introduced by Lutwak \cite{Lu3, LYZ}.
Let $\Pi_p K$, $p \ge 1$ denote the compact convex set whose support function is given by
\begin{equation}\label{eq:pr}
h_{\Pi_p K} (\xi )^p=\frac{1}{2n}\int\limits_{S^{n-1}}|\langle x,\xi\rangle |^p d\,S_p(K,x),
\qquad \xi\in S^{n-1}.
\end{equation}
We note that  $\Pi_1K =n \Pi K$. Moreover, for some fixed $p \ge1 $, we say that $D$ is a $p$-projection body
if $D$ is the $p$-projection body of some convex body. By (\ref{eq:pr}), for every $p$-projection body $D$
there exists a finite measure $\nu_D$ on $S^{n-1}$ such that
\begin{equation}\label{p-proj-body}
h^p_D(x)=\int_{S^{n-1}} |\langle x,\xi\rangle|^p\ d\nu_D(\xi),\qquad x\in {\mathbb R}^n.
\end{equation}
Let us  denote by $\Pi_{p,n}$ the class of all $p$-projection bodies in $\R^n.$

We may pass now to the proof of Theorem~\ref{main-proj}.

\medskip

\noindent {\bf Proof of Theorem~\ref{main-proj}.} Let $\e>0$ be such that for every $\xi\in S^{n-1}$
\begin{equation}\label{proj1}
h^p_{\Pi_p K} (\xi ) \le \e h^p_{\Pi_p L} (\xi ).
\end{equation}
By  (\ref{eq:pr}), the condition (\ref{proj1}) is equivalent to
\begin{equation}\label{proj2}
\int_{S^{n-1}}\!\!\! |\langle x,\xi\rangle |^p\ dS_p(K, x)\le \e \int_{S^{n-1}}\!\!\! |\langle x,\xi\rangle|^p\ dS_p(L,x),\, \forall \xi\in S^{n-1}.
\end{equation}
For small $\delta>0,$ let $D\in \Pi_{p,n}$ be such that $D\subset L$ and
\begin{equation}\label{dist}
|L|^{\frac 1n}\le (1+\delta)d_{\rm vr}(L,\Pi_{p,n})\ |D|^{\frac 1n},
\end{equation}
and let $\nu_D$ be the measure on $S^{n-1}$ corresponding to $D$ by (\ref{p-proj-body}).
Integrating both sides of (\ref{proj2}) with respect to $d\nu_D(\xi),$ we get
$$
\int_{S^{n-1}} \int_{S^{n-1}}\!\!\! |\langle x,\xi\rangle |^p\ dS_p(K, x) d\nu_D(\xi)\le \e \int_{S^{n-1}}\int_{S^{n-1}}\!\!\! |\langle x,\xi\rangle|^p\ dS_p(L,x)d\nu_D(\xi),
$$
for all  $\xi\in S^{n-1}$. Applying Fubini's theorem on $S^{n-1}$ together with (\ref{eq:pr}) we get
\begin{equation}\label{proj3}
\int_{S^{n-1}} h^p_D(x)\ dS_p(K, x)\le\e \int_{S^{n-1}} h^p_D(x)\ dS_p(L, x).
\end{equation}
Since $D\subset L,$ we have $h_D(x)\le h_L(x)$ for every $x\in S^{n-1},$ so the right-hand side of (\ref{proj3})
can be estimated from above by
$$\e \int_{S^{n-1}} h^p_D(x)\ dS_p(L, x)\le \e \int_{S^{n-1}} h^p_L(x)\ dS_p(L, x)= \e n|L|.$$
By (\ref{mixed}), (\ref{pmin}),  (\ref{mink}) and (\ref{dist}), the left-hand side of (\ref{proj3}) can be estimated from below by
\begin{align*}\int_{S^{n-1}} h^p_D(x)\ dS_p(K,x) &=n V_p(K, L) \ge |K|^{\frac{n-p}{n}} |D|^{\frac{p}{n}}  \\
&\ge \frac n{(1+\delta)^pd^p_{\rm vr}(L,\Pi_{p,n})}|K|^{\frac {n-p}n} |L|^{\frac pn}.\end{align*}
Combining these estimates we see that
$$ \frac n{(1+\delta)^pd^p_{\rm vr}(L,\Pi_{p,n})}|K|^{\frac {n-p}n} |L|^{\frac pn}\le\e n|L|.$$
Sending $\delta\to 0,$ we get
$$\left(\frac{|K|}{|L|}\right)^{\frac{n-p}n}\le d^p_{\rm vr} (L, \Pi_{p,n})\ \e.$$
Now, putting
$$\e= \max_{\xi\in S^{n-1}}  \frac{h^p_{\Pi_pK} (\xi)}{h^p_{\Pi_p L} (\xi)},$$
we get the result. $\hfill\Box $

\medskip

A mixed version of the Busemann-Petty and Shephard problems was posed by Milman and solved in \cite{GK1}. Namely,
if $K$ is a convex body in ${\mathbb R}^n$, $D$ is a compact subset of ${\mathbb R}^n$
and $1\le k\le n-1$, then the inequalities
$|K\vert H|\le |D\cap H|$ for all $H\in Gr_{n-k}$ imply $|K|\le |D|.$ Here $K\vert H$ is the orthogonal projection of $K$ onto $H.$
One can easily modify the proof from \cite{GK1} to get a slightly stronger version of this result.

\begin{theorem} Let $K$ be a convex body in $\R^n,$ let $D$ be a compact set in $\R^n,$ and $1\le k \le n-1.$
Then
$$\left(\frac{|K|}{|L|}\right)^{\frac {n-k}n}\le \max_{H\in Gr_{n-k}} \frac{|K\vert H|}{|L\cap H|}.$$
\end{theorem}

\section{Applications}\label{section:5}

\subsection{Comparison theorem for the Radon transform} We can recover the isomorphic Busemann-Petty theorem for the Radon transform established in \cite{KPZv}, as follows. If, in addition to the conditions of Theorem \ref{quotient-main}, we assume that
$$\int_{K\cap H}f\le \int_{L\cap H} g,\qquad \forall H\in Gr_{n-k},$$
then we get
$$\int_K f \le \frac n{n-k} \left(d_{\rm ovr}(K,BP_k^n)\right)^k |K|^{\frac kn}\left(\int_Lg\right)^{\frac{n-k}n}.$$

\subsection{A lower estimate for the sup-norm of the Radon transform} Theorem \ref{quotient-main} with $L=B_2^n$ and $g\equiv 1$ is the slicing inequality for arbitrary functions from \cite{K4} (see inequality (\ref{main-problem1}) above):
$$\int_Kf\le \frac n{n-k} \frac{|B_2^n|^{\frac{n-k}n}}{|B_2^{n-k}|} \left(d_{\rm ovr}(K,BP_k^n)\right)^k |K|^{\frac kn}\max_H \int_{K\cap H} f.$$
Note that the constant $ \frac{|B_2^n|^{\frac{n-k}n}}{|B_2^{n-k}|}$ is less than 1, and $\frac n{n-k}\le e^k.$

\subsection{Mean value inequality for the Radon transform} Let $K=L,$ and $g\equiv 1.$ Then, as we have mentioned in the Introduction,
$$\frac{\int_Kf}{|K|}\le \frac n{n-k} \left(d_{\rm ovr}(K,BP_k^n)\right)^k \max_H \frac{\int_{K\cap H}f}{|K\cap H|}.$$

\subsection{The isomorphic Busemann-Petty problem for sections of proportional dimensions.}
Theorem \ref{quotient-holder} and inequality (\ref{kpzy}) imply the following result from \cite{K6}, which solves the isomorphic Busemann-Petty problem in affirmative for sections of proportional dimensions. If $K,L$ are origin-symmetric convex bodies in
$\R^n$ and $k\ge \lambda n,$ where $0<\lambda <1,$ so that $|K\cap H|\le |L\cap H|$ for every $H\in Gr_{n-k},$ then
$|K|^{\frac{n-k}n}\le (C(\lambda))|L|^{\frac{n-k}n},$ where the constant $C(\lambda)$ depends only on $\lambda.$

\subsection{Inequalities for projections.} Theorem \ref{main-proj} immediately implies an isomorphic version of the Shephard problem first established by
Ball; it immediately follows from \cite{Ba2,Ba3}.

\begin{co} \label{ball} Let $K$ and $L$ be origin-symmetric
convex bodies in $\R^n$ such that
$$|K\vert\xi^\bot|\le |L\vert\xi^\bot|,\qquad \forall \xi\in S^{n-1},$$
then
$$|K|\le d_{\rm vr}(L,\Pi_n) |L|.$$

\end{co}

\begin{co}\label{slicing-proj} Let $L$ be an origin-symmetric convex body in $\R^n.$ Then
$$\min_{\xi\in S^{n-1}} |L\vert\xi^\bot|\le \sqrt{e} \ d_{\rm vr}(L,\Pi_n)\  |L|^{\frac{n-1}n}.$$
\end{co}

\pf Apply Theorem \ref{main-proj} to $K=B_2^n$ and $L,$
and then use the fact that $c_{n,1}=|B_2^{n-1}|/|B_2^n|^{\frac{n-1}n}\le \sqrt{e}.$ \endpf

By John's theorem \cite{J} and the fact that ellipsoids are projection bodies (see, for example \cite{K1, S2}), we have $d_{\rm vr}(L,\Pi_n)\le \sqrt{n}$
for any origin-symmetric convex body $L$ in $\R^n.$ On the other hand, Ball \cite{Ba2} proved that
there exists an absolute constant $c>0$ so that for every $n$ there exists an origin-symmetric convex
body $L_n$ of volume $1$ in $\R^n$ satisfying $|L\vert\xi^\bot|\ge c\sqrt{n}$ for all $\xi\in \R^n.$ Combined with
Corollary \ref{slicing-proj}, these estimates show that
\begin{equation}\label{distance-proj}
c\sqrt{n}\le \max_L d_{\rm vr}(L,\Pi_n)\le \sqrt{n},
\end{equation}
where $c$ is an absolute constant, and the maximum is taken over all origin-symmetric convex bodies in $\R^n.$
This estimate was first established by Ball; it immediately follows from \cite[Example 2]{Ba3}.

Note that the distance $d_{\rm vr}(L,\Pi_n)$ has been studied by several authors. It was introduced in \cite{Ba3} and was proved to be equivalent to the {\it weak-right-hand-Gordon-Lewis} constant of $L$. Also it was connected to the random unconditional constant of the dual space (see Theorem 5 and Proposition 6 in \cite{Ba3}). In \cite{GMP} this distance was called zonoid ratio, and it was proved that it is bounded from above by the projection constant of the space. In the same paper  the zonoid ratio was computed for several classical spaces. We refer the interested reader to \cite{GMP}, \cite{GJN}, \cite{GJ}, \cite{GLSW} for more information.

\subsection{Milman's estimate for the isotropic constant}
We say that a compact set $K$ with volume 1 in $\R^n$ is in isotropic position if for each $\xi\in S^{n-1}$
$$\int_K \langle x,\xi\rangle^2 dx = L_K^2$$
where $L_K$ is a constant that is called the isotropic constant of $K.$ In the case where $K$ is origin-symmetric convex, the slicing problem of Bourgain is equivalent to proving that $L_K$ is bounded by an absolute constant.

Hensley \cite{He} has proved that there exist absolute constants $c_1,c_2>0$ so that for any origin-symmetric convex body $K$ in $\R^n$ in isotropic position and any $\xi\in S^{n-1}$
$$\frac{c_1}{L_K}\le |K\cap \xi^\bot| \le \frac{c_2}{L_K}.$$

The following inequality was proved by Milman \cite{M1}. We present a simpler proof.
\begin{theorem} There exists an absolute constant $C$ so that for any origin-symmetric isotropic convex body $K$ in $\R^n$
$$L_K\le C\  d_{\rm ovr}(K,{\cal{I}}_n).$$
\end{theorem}

\pf By Theorem \ref{quotient-holder} with $k=1$ and Hensley's theorem, for any origin-symmetric isotropic convex bodies $K,D$
in $\R^n$
$$\left(\frac{|K|}{|D|}\right)^{\frac {n-1}n}\le d_{\rm ovr}(K,{\cal{I}}_n) \max_{\xi\in S^{n-1}} \frac{|K\cap \xi^\bot|}{|D\cap \xi^\bot|} \le d_{\rm ovr}(K,{\cal{I}}_n) \frac{|K|^{\frac {n-1}n} \frac{c_2}{L_K}} {|D|^{\frac {n-1}n} \frac{c_1}{L_D}},$$
where $c_1,c_2>0$ are absolute constants, so
$$\frac{L_K}{L_D} \le C\ d_{\rm ovr}(K,{\cal{I}}_n).$$
Now put $D=B_2^n / |B_2^n|^{\frac 1n},$  and use the fact that $L_D$ is bounded by an absolute constant;
see \cite{BGVV}. \endpf

\end{document}